\documentclass[12pt]{article}
\usepackage{amsfonts}
\usepackage{amsthm}
\usepackage{graphicx}
\usepackage{latexsym}
\makeatletter
\renewcommand{\@seccntformat}[1]
{\csname the#1\endcsname.\enspace}
\makeatother
\renewcommand{\baselinestretch}{1}
\begin{document}
\newtheorem{theorem}{Theorem}
\newtheorem{lemma}{Lemma}
\newtheorem{remark}{Remark}
\newtheorem{corollary}{Corollary}
\newtheorem{example}{Example}
\renewcommand{\theequation}{\thesection.\arabic{equation}}
\renewcommand{\thetheorem}{\thesection.\arabic{theorem}}
\renewcommand{\thelemma}{\thesection.\arabic{lemma}}
\renewcommand{\thecorollary}{\thesection.\arabic{corollary}}
\renewcommand{\theremark}{\thesection.\arabic{remark}}
\renewcommand{\theexample}{\thesection.\arabic{example}}
\begin{center}
\textbf{Estimation of a multivariate normal mean with a bounded signal to noise ratio} 
%\footnote{draft of \today}
\\
\end{center}

\begin{center}
{Othmane Kortbi, \' Eric
 Marchand\footnote{Corresponding author: eric.marchand@usherbrooke.ca} } \\
\end{center}
\begin{center}
{\it  Universit\'e de
    Sherbrooke, D\'epartement de math\'ematiques, Sherbrooke, QC,
    CANADA, J1K 2R1} \\
\end{center}

\begin{center}
{\sc Summary}
\end{center}
\small  
For normal canonical models with $X \sim N_p(\theta,  \sigma^{2} I_{p}), \;\; S^{2} \sim \sigma^{2}\chi^{2}_{k}, \;\mbox{independent}$, we consider the problem of estimating $\theta$ under scale invariant squared error loss $\frac{\|d-\theta \|^{2}}{\sigma^{2}}$, when it is known that the signal-to-noise ratio $\frac{\left\|\theta\right\|}{\sigma}$ is bounded above by $m$.  Risk analysis is achieved by making use of a conditional risk decomposition and we obtain in particular sufficient conditions for an estimator to dominate either the unbiased estimator $\delta_{UB}(X)=X$, or the maximum likelihood estimator $\delta_{\hbox{mle}}(X,S^2)$, or both of these benchmark procedures.  The given developments bring into play the pivotal role of the boundary Bayes estimator $\delta_{BU}$ associated with a prior on $(\theta,\sigma)$ such that $\theta|\sigma$ is uniformly distributed on the (boundary) sphere of radius $m$ and a non-informative $\frac{1}{\sigma}$ prior measure is placed marginally on $\sigma$.  With a series of technical results related to $\delta_{BU}$; which relate to particular ratios of confluent hypergeometric functions; we show that, whenever $m \leq \sqrt{p}$ and $p \geq 2$, $\delta_{BU}$ dominates both $\delta_{UB}$ and $\delta_{\hbox{mle}}$.  The finding can be viewed as both a multivariate extension of $p=1$ result due to Kubokawa (2005) and a unknown variance extension of a similar dominance finding due to Marchand and Perron (2001).  Various other dominance results are obtained, illustrations are provided and commented upon.  In particular, for $m \leq \sqrt{\frac{p}{2}}$, a wide class of Bayes estimators, which include priors where $\theta|\sigma$ is uniformly distributed on the ball of radius $m$, are shown to 
dominate $\delta_{UB}$.


\normalsize

\section{Introduction}
\setcounter{equation}{0} 
\setcounter{theorem}{0}
\setcounter{lemma}{0} 
\setcounter{corollary}{0}
\setcounter{remark}{0}
\subsection{The model}
We consider the normal canonical model  
\begin{equation}\label{model}
X \sim N_p(\theta,  \sigma^{2} I_{p}), \;\; S^{2} \sim \sigma^{2}\chi^{2}_{k}, \;\mbox{independent} \,,
\end{equation}
with $p \geq 1$ and $k \geq 1$,  
which plays a central role in both statistical theory and practice.  We wish to estimate $\theta$ under scale invariant squared error loss
\begin{equation}
\label{loss}
L\left(\left(\theta, \sigma\right), d\right)=\frac{\|d-\theta \|^{2}}{\sigma^{2}} ,
\end{equation}
and where $\left(\theta, \sigma\right)$ is restricted to the parameter space 
\begin{equation}\label{ps}
\Theta\left(m\right)=\left\{\left(\theta, \sigma\right) \in \Re^{p}\times \Re^{+}: \frac{\left\|\theta\right\|}{\sigma}\leq m \,, 0<\sigma< \infty\right\}\,,
\end{equation}
where $m$ is some known positive constant.
Viewing $\frac{\left\|\theta\right\|}{\sigma}$ as a multivariate version of a signal to noise ratio, the problem can be described as estimating a normal mean with upper bounded signal to noise ratio.  Alternatively, as previously described and
analyzed by Kubokawa (2005) for $p=1$, the parametric constraint places an upper bound on the reciprocal of the
coefficient of variation in absolute value.  

We will be concerned with Bayesian inference in such restricted
parameter space problems, which does not, conceptually, present any difficulties as both the prior and the resulting posterior (it it exists) will be
adapted and will adapt to the constraints.  Assessing the frequentist performance of Bayesian estimators in such situations is, however,
considerably more challenging.  Such assessments may include, for instance, testing for minimaxity, an evaluation in comparison
to a benchmark procedure such as minimum risk equivariant (MRE) estimator or a maximum likelihood estimator (mle), or a study of the frequentist performance
of associated Bayesian credible sets.

Model (\ref{model}) and the constraint (\ref{ps}) arise in normal full rank linear models $Y \sim N_n(Z\beta, \sigma^2 I_n)$
with orthogonal design matrix $Z (n \times p)$; $n > p$; unknown parameter vector $\beta (p \times 1)$, and with the constraint
$\frac{\left\|\beta\right\|}{\sigma} \leq m$.  Indeed, the sufficient statistic $(\hat{\beta}=(Z'Z)^{-1} Z'Y = Z'Y, W=\|Y - Z \hat{\beta} \|^2)$ has a distribution which matches (\ref{model}) with $X= \hat{\beta}$, $S^2=W$, $k=n-p$ accompanied by constraint (\ref{ps}) with $\theta=\beta$.  More generally, where $Z$ is not necessarily orthogonal, the correspondence also arises as above through the constraint $\frac{\|Z\beta \|}{\sigma} \leq m$ by setting $X=(Z'Z)^{-1/2} \hat{\beta}$ and 
$\theta= (Z'Z)^{-1/2} \beta$.
 
As considered and motivated by Kariya, Giri and Perron (1988), as well as Perron and Giri (1990), model
(\ref{model}) arises for the curved model setting
\begin{equation}
\label{kgpmodel}
Y_1, \ldots, Y_n \; \mbox{independent} \; N_p(\mu, \sigma^2 I_p)\,, \; \hbox{with} \; \sigma^2= \frac{\mu' \mu}{\lambda}.
\end{equation}  
The correspondence is achieved by considering the sufficient statistic $\bar{Y}= \frac{1}{n} \sum_{i=1}^n \, Y_i$, 
$W= tr (\sum_{i=1}^n \, (Y_i- \bar{Y})(Y_i- \bar{Y})')$ and setting $X=\sqrt{n} \, \bar{Y}$, $S^2=W$, $\theta= \sqrt{n} \mu$,
and $k=(n-1)p$ in (\ref{model}).  Constraint (\ref{ps}) arises by assuming $\lambda \in [0,m]$ in (\ref{kgpmodel}), while the findings of 
Kariya, Giri and Perron (1988), as well as Perron and Giri (1990), relate to known $\lambda$ and the best equivariant estimator
$\delta_{\lambda}(X,S)$ for loss $\frac{\|d - \mu \|^2}{\mu' \mu}$ (equivalent to \ref{loss}) under the group of transformations $G= \Re^+ \times O(p)$, $\Re^+$ being the multiplicative group of positive real numbers and $O(p)$ being the group of $p \times p$ orthogonal matrices.  We will extract some features of the invariance structure and of the estimator 
$\delta_{\lambda}(X,S)$ in Subsection 2.2, but refer to the above mentioned papers for additional details.

\subsection{The estimators and previous findings}
 
Two benchmark procedures are the unbiased estimator $\delta_{UB}(X, S^2)= X$ and the maximum likelihood estimator $\delta_{\hbox{mle}}(X, S^2)$.  Although both of the above procedures will be shown to be inadmissible (for $p \geq 1$ 
and $p \geq 2$ respectively), it is nevertheless of interest to describe the structure of improvements with an emphasis on potential Bayesian improvements. 

Actually, the choice $\delta_{UB}(X, S^2)$ is clearly inefficient for our problem, as 
seen for instance by a straightforward analysis of the risk of affine linear estimators of the form $\delta_a(X,S^2)=aX; a \in \Re$.  Indeed, the choices $\delta_a$ with $a \in [\frac{m^2-p}{m^2+p}, 1)$ dominate $\delta_{UB}(X, S^2)$ and the minimax procedure in this class, given by the choice $a=\frac{m^2}{m^2+p}$, always dominates (strictly) $\delta_{UB}(X, S^2)$.  
We will show below that $\delta_{\hbox{mle}}(X, S^2)$, which takes into account the constraint in opposition to $\delta_{UB}(X, S^2)$, also dominates $\delta_{UB}(X, S^2)$ for all $(m,p)$, so that it is more challenging to obtain improvements on the former as opposed to the latter.

\begin{remark} (Minimaxity)
The above paragraph also implies that $\delta_{UB}(X, S^2)$ is never minimax for our problem, in contrast with the case where $\theta$, rather than $\frac{\theta}{\sigma}$, is constrained to a ball of radius $m$ and the loss is the same as in (\ref{loss}). Indeed, as shown by Kubokawa (2005) for $p=1$ and as a consequence of a general minimax result given by Marchand and Strawderman (2011) applicable for $p \geq 1$, $\delta_{UB}(X, S^2)$ is minimax (but still inadmissible) for the restricted parameter space with $\|\theta\| \leq m$ and $\, 0<\sigma< \infty$ .
\end{remark}

We denote $S_h$ and $B_h$, respectively, as the $p$ dimensional sphere and ball of radiuses $h$ centered at the origin.
We consider priors on the restricted parameter space $\Theta(m)$ admitting the representations:
\begin{equation}
\label{prior}
{\bf (i)} \; \theta|\sigma^2 \sim \hbox{Uniform}(S_{m \sigma}), \sigma^2 \sim (\sigma^2)^{\frac{l}{2}-1}, \; 
\hbox{or} \; {\bf (ii)}  \; \theta|\sigma^2 \sim \hbox{Uniform}(B_{m \sigma}), \sigma^2 \sim (\sigma^2)^{\frac{l}{2}-1}\,,
\end{equation}
with $l < p+k\,$ so that the associated posterior distributions be well defined.  These interesting priors are improper, but they are proper with respect to $\theta$ for fixed $\sigma$.  The case $l=0$ is of particular significance since it corresponds to the right Haar invariant case for $\sigma^2$, and since the prior in {\bf (ii)} corresponds to the truncation onto the restricted parameter space $\Theta(m)$ of the usual noninformative prior for $(\theta, \sigma^2)$ set on $\Re^{p}\times \Re^{+}$.  The Bayes estimators associated with priors in {\bf (i)} and {\bf (ii)} may hence be described as boundary uniform or fully uniform, and we denote these $\delta_{BU,\,l}$ and $\delta_{U,l}$ respectively.

The above priors and Bayes estimators are analogous to the estimators $\delta_{BU}$ and $\delta_{U}$ studied by Marchand and Perron (2001) for the known $\sigma^2$ case in model (\ref{model}), with loss (\ref{loss}) and the restriction (\ref{ps}), where both $\delta_{BU}$ and $\delta_{U}$ were proven to dominate the maximum likelihood estimator $\delta_{\hbox{mle}}(X) =
(\frac{m}{\|X\|} \wedge 1) X$ for sufficiently small $m$; and where namely $\delta_{BU}$ dominates $\delta_{\hbox{mle}}(X)$, and hence $X$, under the simple condition $m \leq \sqrt{p}$.  Furthermore, a result by Hartigan (2004) applicable for much more general restricted convex parameter spaces, and reviewed along with analogous findings by Marchand and Strawderman (2004), implies that $\delta_{U}$ dominates $X$ for all $(m,p)$.  

Returning to our unknown $\sigma^2$ case, Kubokawa (2005) provided various improvements on $\delta_{UB}(X, S^2)=X$ for 
the univariate case $p=1$, showing {\bf (a)} dominance of $\delta_{U,l}$ for all $m>0$ and $l \in [0, k+1]$, as well as {\bf (b)} the dominance of $\delta_{BU,\,l}$ for $m \leq 1$ and $l \in [0, k+1]$.  Hence, these findings may be viewed as 
unknown $\sigma^2$ univariate extensions of some of the existing results in the literature, including some described in the 
previous paragraph.  However, Kukokawa's results are limited to the univariate case and do not apply to either the maximum likelihood estimator $\delta_{\hbox{mle}}(X, S^2)$, nor to affine linear estimators $aX$, $a \in [0,1)$.

In this paper, we extend Kubokawa's findings concerning $\delta_{BU,\,l}$ to the multivariate case $p>1$ providing sufficient conditions for dominance of both $\delta_{UB}(X, S^2)$ and $\delta_{\hbox{mle}}(X, S^2)$.  Namely, we show that: {\bf (i)} the Bayes estimators $\delta_{BU,\,l}$, $ 0 \leq l < k+p$, dominate $\delta_{UB}(X, S^2)$ whenever 
$m \leq \sqrt{p}$, and {\bf (ii)} $\delta_{BU,\,0}$ dominates $\delta_{\hbox{mle}}(X, S^2)$ whenever $m \leq \sqrt{p}$ and $p \geq 2$ (Corollary \ref{riskcomparison}).  This yields a striking parallel with Marchand and Perron's dominance result for the known $\sigma^2$ case.  We also show that $\delta_{\hbox{mle}}$ always improves on $\delta_{UB}$, infer other dominance results, and present illustrations and accompanying commentary.  For very small parameter spaces (precisely $m \leq \sqrt{\frac{p}{2}}$), we obtain a universal dominance result showing that a vast class of Bayes estimators, which includes $\delta_{U,l}$ and $\delta_{BU,l}$ for $l \leq 0$, dominate $\delta_{UB}$.

Our methods also depart from, and arguably simplify, the methods used by Kubokawa in the univariate case applicable to $\delta_{BU,l}$.   Key features of Kubokawa's results are techniques previously introduced by Kubokawa himself (Kubokawa, 2004, Integral Expression for Risk Difference (IERD), as well as Marchand and Strawderman (2005).  With extensions of these techniques 
challenging to obtain and elusive in the multivariate case, with the absence of results applicable to other estimators such as
the $\delta_{\hbox{mle}}(X, S^2)$, we rather exploit a conditional risk decomposition on a maximal invariant $T= \frac{\|X\|^2}{S^2}$ in a similar fashion as Marchand and Perron (2001) and Moors (1985).  The deficiency of $\delta_{\hbox{mle}}(X, S^2)$ or $\delta_{\hbox{mle}}(X, S^2)$ is revealed as one of providing estimates that are too far from the origin.  Even when $m > \sqrt{p}$, we will show that improvements necessarily occur by projecting towards the Bayes estimator $\delta_{BU,0}$.
The Bayes estimator $\delta_{BU,0}$ coincides with the best equivariant estimator $\delta_{m}$, and 
analytical results for $\delta_{BU,0}$, which we will describe and make use of, were 
previously given by Kariya, Giri and Perron (1988), as well as Perron and Giri (1990).  

The remainder of this paper is organized as follows.   Section 2 contains key features and properties of the invariance structure and the risk function of equivariant estimators, as well as various expressions and key properties relative to the Bayes estimators $\delta_{BU,\,l}$ and the maximum likelihood estimator $\delta_{\hbox{mle}}(X, S^2)$.  The dominance results are presented in Section 3 along with various illustrations and comments.  Final remarks are given in Section 4.

\section{Preliminary results}

\subsection{The invariance structure and risk function}

The challenge here is to obtain good improvements or alternatives under loss (\ref{loss}) that 
capitalize on the parametric information in (\ref{ps}), as measured by the risk function 
\begin{equation}\label{risk}
	R((\theta, \sigma), \delta)=\frac{1}{\sigma^{2}} E_{\theta, \sigma}\left[\left\|\delta(X, S)-\theta\right\|^{2} \right] 
	,  \left(\theta, \sigma\right) \in  \Theta\left(m\right)
	\,.
\end{equation}
We provide findings for (nonrandomized) equivariant estimators, with respect to the group structure described in Section 1.1, and as shown by Kariya, Giri and Perron (1988), as well as Perron and Giri (1990), to be of the form
\begin{equation}
\label{deltah}
\delta_h(X, S^2) = h(\frac{\|X\|^2}{S^2} ) X\,,
\end{equation}
for some measurable function $h(\cdot)$.  Equivariant estimators are thus collinear to $X$ and conveniently represented
by the corresponding multiplier $h$ which controls the degree of expansion or shrinkage with respect to $X$, depending on $(X,S^2)$ only through the maximal invariant statistics $T= \frac{\|X\|^2}{S^2}$.   The class of such estimators include $\delta_{UB}$ ($h \equiv 1$), affine linear estimators $aX$ ($h \equiv a$), $\delta_{\hbox{mle}}$ (see Lemma \ref{mle}),  generalized Bayes estimators $\delta_{BU,\,l}$ (see Lemma \ref{bul}), and more generally Bayes estimators with a spherically symmetric structure which includes $\delta_{U,\,l}$ (see Lemma \ref{ss}).
Equivariant estimators have a risk function in (\ref{risk}) depending
 on the unknown parameters only through the maximal invariant $\lambda= \frac{\|\theta\|}{\sigma} \in [0,m]$.  By a slight abuse of notation, we will write this risk
$R(\lambda, \delta_h)$.  Here is a useful representation for the risk $R(\lambda, \delta_h)$ 
of equivariant estimators achieved by conditioning on the maximal invariant statistic $T$ and highlighting the key role of the best equivariant estimator $\delta_{\lambda}(X,S^2)$ (for $\frac{\|\theta\|}{\sigma}=\lambda$).

\begin{lemma}  
\label{risklemma}
Under model (\ref{model}) and loss (\ref{loss}), we have
\begin{equation}
\label{riskdecomposition}
	R(\lambda, \delta_{h}) =  
	\lambda^{2}+E_{\lambda} \left[ a(T) \left\{\left(h(T)-h_{\lambda}(T)\right)^{2}-h_{\lambda}^{2}(T)\right\} \right]
	\, ,
\end{equation}
with $a(T)=E_{\theta, \sigma}\left( \left. \frac{X' X}{\sigma^{2}} \right| T \right)$,
$b(T)= E_{\theta, \sigma} \left( \left. \frac{\theta' X}{\sigma^{2}} \right| T   \right)$, and
$h_{\lambda}(T)= \frac{b(T)}{a(T)}$.
\end{lemma}
{\bf Proof.}  The result is immediate by writing the loss as $\frac{\|h(t) x - \theta \|^2}{\sigma^2} = h^2(t) \frac{x'x}{\sigma} + \lambda^2 - 2 h(t) \frac{\theta'x}{\sigma^2 }$ and decomposing the risk as $E[L((\theta, \sigma^2), \delta_h)]
= E^T\{E[L((\theta, \sigma^2), \delta_h)]|T \}$.  \qed

\begin{remark}
As seen by (\ref{riskdecomposition}) above, the risk of $\delta_h$ is constant for the restriction $\frac{\|\theta\|}{\sigma}=\lambda$ and the optimal procedure is given by the BEE $\delta_{\lambda}(X,S^2)=h_{\lambda}(T) X$.
\end{remark}

Kariya, Giri and Perron (1988), as well as Perron and Giri (1990) gave an explicit expression for the BEE $\delta_{\lambda}(X,S^2)$.  It is reproduced below in
Lemma \ref{bul} where we derive an expression for the Bayes estimators $\delta_{BU,\,l}$.  The estimators $\delta_{m}(X,S^2)$ and $\delta_{BU,\,l}$ necessarily coincide by virtue of general relationships between best equivariant estimators and Bayes estimators with respect to Haar right invariant priors (e.g., Eaton, 1989).

%\begin{lemma}\label{T2} \textbf{(Kariya, Giri et Perron, 1988)}.
%For $\frac{\left\|\theta\right\|}{\sigma}=\lambda$, the best equivariant estimator, 
%$\delta_{\lambda}(x, s)=h_{\lambda}\left( t \right) x$, is determined by
%\begin{equation} \label{4}
% h_{\lambda}(t)= \frac{\lambda^{2}}{p} \frac{F\left(\frac{k+p}{2}+1, \frac{p}{2}+1, \frac{\lambda^{2}t}{2(1+t)}\right)}
% {F\left(\frac{k+p}{2}+1, \frac{p}{2}, \frac{\lambda^{2}t}{2(1+t)}\right)}
%\, ,
%\end{equation}
%\end{lemma}

\subsection{Bayes estimators}

We begin here with a general expression for Bayes estimators associated with priors of the form
\begin{equation}\label{generalprior}
\pi\left(\theta, \sigma^{2}\right) =\pi_{\sigma}\left(\theta\right) \sigma^{l-2} \,, \;\;\;\; l < k+p \,,
\end{equation} 
where $\pi_{\sigma}(\cdot)$ is for fixed $\sigma$ a (proper) density with respect to a finite measure $\nu_{\sigma}$
supported on, or a subset of, the ball $B(m\sigma)$.

\begin{lemma}\label{bayeslemma}
For model (\ref{model}) and loss (\ref{loss}), 
Bayes estimators with respect to priors as in (\ref{generalprior}) are equal to
\begin{equation}\label{deltapi}
\delta_{\pi}\left(x, s^2\right)= x+\frac{\nabla_{x} \int_{0}^{\infty}\sigma^{l-k-p-2} \ e^{-\frac{s^{2}}{2\sigma^{2}}}  \, m^{\pi}_{\sigma}\left(x\right) \,d\sigma^{2}}
{\int_{0}^{\infty} \sigma^{l-k-p-4} \ e^{-\frac{s^{2}}{2\sigma^{2}}} \, m^{\pi}_{\sigma}\left(x\right) \,d\sigma^{2}}\, ,
\end{equation}
where $\nabla_{x}$ denotes the gradient vector with respect to $x$, and $(2\pi \sigma^2)^{-p/2} \, m^{\pi}_{\sigma}$
is the marginal density of $X|\sigma$ with
\begin{equation}\label{mpisigma}
m^{\pi}_{\sigma}\left(x\right)=\int_{B_{m\sigma}} e^{-\frac{1}{2\sigma^{2}}\left\|x-\theta\right\|^{2}}  \pi_{\sigma}\left(\theta\right) \,d\nu_{\sigma}(\theta) \, .
\end{equation}
\end{lemma}
\noindent {\bf Proof.}
We have
\begin{eqnarray}
\label{formula} \delta_{\pi}\left(x, s^2\right) &=& \frac{E\left[\ \theta/\sigma^{2}\left| x, s^2 \right.\right]}
{E\left[\ 1/\sigma^{2}\left| x, s^2 \right.\right]}
\\
\label{ar} &=& \frac{\int_{0}^{\infty}\int_{B_{m\sigma}}\theta \ e^{-\frac{1}{2\sigma^{2}}\left\|x-\theta\right\|^{2}}e^{-\frac{s^{2}}{2\sigma^{2}}}\ \pi_{\sigma}\left(\theta\right) \sigma^{l-k-p-4} \,d\nu_{\sigma}(\theta) \,d\sigma^{2}}
{\int_{0}^{\infty}\int_{B_{m\sigma}} e^{-\frac{1}{2\sigma^{2}}\left\|x-\theta\right\|^{2}}e^{-\frac{s^{2}}{2\sigma^{2}}}\ \pi_{\sigma}\left(\theta\right) \sigma^{l-k-p-4} \,d\nu_{\sigma}(\theta) \,d\sigma^{2}}
\\
\nonumber &=& x+
 \frac{\int_{0}^{\infty}\sigma^{l-k-p-4} \ e^{-\frac{s^{2}}{2\sigma^{2}}} \int_{B_{m\sigma}}\left(\theta -x\right)\ e^{-\frac{1}{2\sigma^{2}}\left\|x-\theta\right\|^{2}} \pi_{\sigma}\left(\theta\right) \,d\nu_{\sigma}(\theta) \,d\sigma^{2}}
{\int_{0}^{\infty} \sigma^{l-k-p-4} \ e^{-\frac{s^{2}}{2\sigma^{2}}} \int_{B_{m\sigma}} e^{-\frac{1}{2\sigma^{2}}\left\|x-\theta\right\|^{2}}  \pi_{\sigma}\left(\theta\right) \,d\nu_{\sigma}(\theta) \,d\sigma^{2}} \, .
\end{eqnarray}
The result follows since
\begin{equation}
\nonumber
\label{82}
\sigma^2  \nabla_{x} m^{\pi}_{\sigma}\left(x\right)=\int_{B_{m\sigma}} \left(\theta -x\right)\ e^{-\frac{1}{2\sigma^{2}}\left\|x-\theta\right\|^{2}} \pi_{\sigma}\left(\theta\right) \,d\nu_{\sigma}(\theta) \,.\;\; \;\;\;\;\;\;\;\;\;\; \qed
\end{equation}

Lemma \ref{bayeslemma} applies to the boundary uniform and fully uniform priors in (\ref{prior}), among others.  We pursue here 
with useful expressions for the former, which will also serve as a benchmark for other Bayesian estimators (see Lemma \ref{ss}) in Section 3.

\begin{theorem}\label{bul} We have $\delta_{BU,\,l}\left(X, S^2\right)=h(m,l,T) \, X$ with $T=\frac{\left\|X \right\|^{2}}{S^{2}}$ and
\begin{equation} \label{hbul}
 h(m,l,t) = \frac{m^{2}}{p} \frac{F\left(\frac{k+p-l}{2}+1, \frac{p}{2}+1, \frac{m^{2}t}{2(1+t)}\right)}
 {F\left(\frac{k+p-l}{2}+1, \frac{p}{2}, \frac{m^{2}t}{2(1+t)}\right)} 
\, ,
\end{equation}
$F$ being the confluent hypergeometric function given by $F(a, b, z)=\sum_{i=0}^{\infty}\frac{(a)_{i}}{(b)_{i}}\frac{z^{i}}{i!}$, $z \in \Re$, with $(c)_{i}=\prod_{j=0}^{i-1}(c+j)$ for $i=1, 2, \ldots$, and $(c)_{0}=1$.
\end{theorem}

Here is a familiar identity (e.g., Watson, 1983), related to the normalization constant of a Langevin distribution and useful in the proof of Theorem \ref{bul}.  

\begin{lemma} \label{T12}
For $y \in \Re^p$, and $U$ uniformly distributed on the sphere of radius $r$, we have
\begin{equation}
\label{langevin}
E_r[e^{y'U}] = \Gamma(\frac{p}{2}) \, \frac{I_{p/2-1}\left( \left\|y\right\|r\right)}
{\;\;\;\;\left( \left\|y\right\|r\right)^{(p/2 -1)}} \, ,
\end{equation}
%where $\mathcal U_{r}$ is the uniform measure on the sphere $S_{r}$ of radius $r$, 
where $I_{\nu}\left(\cdot\right)$ is the modified Bessel function of order $\nu$ given by
$$ I_{\nu}(z)= \, \sum_{k \geq 0} \frac{(\frac{z}{2})^{\nu + 2k}}{k! \, \Gamma(\nu + k + 1)}\,.$$
\end{lemma}
\noindent {\bf Proof of Theorem \ref{bul}.}
It suffices to calculate $m^{BU}_{\sigma}\left(x\right) = m^{\pi}_{\sigma}\left(x\right)$ for a boundary uniform prior 
$\pi$ as in (\ref{prior}, {\bf i}).  We have from (\ref{mpisigma}) and (\ref{langevin}), with $U$ uniformly distributed on 
$S_{m}$,
\begin{eqnarray} \label{100}
m^{BU}_{\sigma}\left(x\right)&=& E_m\left[e^{-\frac{1}{2}\|\frac{x}{\sigma} - U\|^2} \right] \nonumber
\\
&=& e^{-\frac{m^{2}}{2}} e^{-\frac{\left\|x\right\|^{2}}{2\sigma^{2}}} E_m\left[ e^{\frac{x'}{\sigma} U} \right] \nonumber
\\
%&=& k \ m^{p-1} \ e^{-\frac{m^{2}}{2}} e^{-\frac{\left\|X\right\|^{2}}{2\sigma^{2}}} \int_{S_{m\sigma}} e^{\frac{X'}{\sigma} %\frac{\theta}{\sigma}}  \,d\mathcal U_{r}\left(\theta\right) \nonumber
%\\
\label{mbu}
&=& \ \ e^{-\frac{m^{2}}{2}} e^{-\frac{\left\|x\right\|^{2}}{2\sigma^{2}}} \, \Gamma(\frac{p}{2}) \,  
{(\frac{m \left\| x\right\|}{\sigma})}^{(1-p/2)}   I_{p/2-1} \left( m\left\|x\right\|\right/\sigma) \, .
\end{eqnarray}
Now, using the identity $\frac{\,d}{\,dt} \left( t^{1-\nu}\, I_{\nu}\left(t\right) \right)= 
t^{1-\nu}\, I_{\nu+1}\left(t\right) $,   calculations lead to 
 \begin{equation}\label{101}
\nabla_{x} m^{BU}_{\sigma}\left(x\right)=- \frac{x}{\sigma^{2}}m^{BU}_{\sigma}\left(x\right)+
\\
\Gamma(\frac{p}{2}) \,  e^{-\frac{m^{2}}{2}} e^{-\frac{\left\|x\right\|^{2}}{2\sigma^{2}}} \, m^{3-p/2} \, \sigma^{p/2-2} \frac{x}{\left\|x\right\|^{p/2-1}}\    I_{p/2}\left( m\left\|x\right\|\right/\sigma) \, .
\end{equation}
Substituting this and (\ref{mbu}) into (\ref{deltapi}) with an interchanging of $\nabla_x$ and $\int$ (justified) yields
\begin{equation}
\delta_{BU,l}\left(x, s^2\right)= m \frac{x}{\left\|x\right\|}\frac{\int_{0}^{\infty}\sigma^{l-k-p/2-4} \  e^{-\frac{\left\|x\right\|^{2}+s^{2}}{2\sigma^{2}}}   I_{p/2}\left( m\left\|x\right\|\right/\sigma) \,d\sigma^{2}}
{ \int_{0}^{\infty}\sigma^{l-k-p/2-5} \  e^{-\frac{\left\|x\right\|^{2}+s^{2}}{2\sigma^{2}}}   I_{p/2-1}\left( m\left\|x\right\|\right/\sigma) \,d\sigma^{2} } \, .
\end{equation}
Finally, the result follows by substituting the series expression for the Bessel functions above, interchanging sums and integrals, integrating out with respect to $\sigma^2$, and some simplification. \footnote{Equivalently, one can proceed with an intermediate identity of the form 
$$ \int_0^{\infty} A^{-\alpha} e^{-T/2A} I_{\nu}(\frac{\mu}{\sqrt{A}})\, dA \,=\, 
\frac{\Gamma(\alpha+\nu/2-1)}{\Gamma(\nu+1)} (\frac{\mu^2}{2T})^{\nu/2} (\frac{2}{T})^{\alpha-1} \,
F(\alpha+\nu/2-1, \nu+1, \frac{\mu^2}{2T} )\,,$$
for positive $\mu, \alpha, T$.}   \qed

%Finally, the result follows by making use of the following identity (e.g., Abramowitz and Stegun, 1964)
%\begin{equation}\label{103}
%F\left(a, c, \frac{u}{z}\right)=\frac{\Gamma\left(c\right)}{\Gamma\left(a\right)} \ z^{a} u^{\frac{1-c}{2}}
%\int_{0}^{\infty} y^{a-\frac{1}{2}-\frac{c}{2}}e^{-zy} I_{c-1}\left(2\sqrt{uy}\right) \,dy \, .  \;\;\;\;\;\; \qed
%\end{equation}

Ratios of confluent hypergeometric functions and their properties hence play an important role
here, as witnessed by Theorem \ref{bul}'s representation of $h(m,l,t)$.  Furthermore, key analytical properties of 
$\delta_{BU,\,l}$ (Lemmas \ref{propertiesh} and \ref{tough}), of other Bayesian estimators (Lemma \ref{ss}) and corresponding risk function comparisons (Section 3)  will hinge on various properties of such ratios, as those given by the following intermediate result.

\begin{lemma}
\label{propertiesratio}
For all $a>0$, $b>0$, $z > 0$, and $c \in \{0,1\}$,  
\begin{enumerate}
\item[{\bf (a)}]  the function  
$ K_{a, b, c}(\cdot)=\frac{F(a-c+1, b-c+1, \cdot)}
 {F(a+1, b, \cdot)}$
is strictly decreasing on $[0,\infty)$ with $\lim_{z \to \infty}K_{a, b, c}(z)= 0$;
\item[{\bf (b)}] the function $ H\left(\cdot\right)=\frac{F\left(\cdot, b+1, z\right)}
 {F\left(\cdot, b, z\right)}$ is strictly decreasing on $(0, \infty)$.
% , and bounded above by $\frac{b}{b+z}$
%on $[b+1, \infty)$.  
\end{enumerate}
\end{lemma}
 \noindent {\bf Proof.}
{\bf(a)} We have for $z \geq 0$ 
$$K_{a, b, c}\left(z\right)=\left(ca+(1-c)b\right) E_{z}\left[\frac{1}{\left(ca+(1-c)b\right)+I}\right] \,,$$ 
where $I$ is a discrete random variable with density proportional to  
$\frac{(a+1)_{i}}{(b)_{i}}\frac{z^{i}}{i!} $ $1_{\left\{0, 1, \ldots\right\}}(i)$.
Since these densities form a family with strictly increasing monotone likelihood ratio in $I$ with parameter $z$, the
result follows since $(ca+(1-c)b > 0$, and $\frac{1}{\left(ca+(1-c)b\right)+i}$ decreases in $i$, $i\in \left\{0, 1, \ldots\right\}$.  The limiting value is established by exploiting the representation (e.g., Abramowitz and Stegun, 1964) 
$F(\alpha, \beta, z)=\frac{\Gamma(\beta)}{\Gamma(\alpha)} \exp(z) z^{\alpha-\beta} \left(1+O\left(\frac{1}{z}\right)\right)$,
where $O\left(\frac{1}{z}\right)$ is a bounded function  by $\frac{1}{z}$ in a neighborhood of infinity. 

{\bf(b)} Similarly, write $\frac{1}{H\left(a\right)}= 1 + \frac{1}{b} \, E_{a}\left[J \right] \,,$
where $J$ is a discrete random variable with mass function 
$$p_{a}(j)\propto \frac{(a)_{j}}{(b+1)_{j}}\frac{z^{j} \exp\{-z\}}{j!} 1_{\left\{0, 1, \ldots\right\}}(j)\,,$$
and observe a strictly increasing monotone likelihood ratio property in $J$ with parameter $a$. This implies the desired monotonicity of $H(\cdot)$.  \qed
% with $$H\left(a\right)\leq H\left(b+1\right)=\frac{b}{b+E_{b+1}\left[J\right]}=\frac{b}{b+z} \,
%\;\,\hbox{for}\;\,a \geq b+1\,.\;\;\;\;\;\;\;\; \qed$$

By making use of these above properties, as well as further properties of confluent hypergeometric functions,
we derive the following analytical results concerning Bayes estimators $\delta_{BU,l}$ given in Theorem \ref{bul}.
More precisely, we now describe how $h(\lambda, l, t)$ varies with respect to $(\lambda,l,t)$.

\begin{lemma}
\label{propertiesh}
\begin{enumerate}
 \item[{\bf (a)}] The function $h(\lambda, l, \cdot)$ is, for $\lambda > 0$, $l < k+p$, strictly decreasing on 
 $\left[0, \infty\right)$ with $\sup_{t\geq 0}h(\lambda, l, t)=h(\lambda, l, 0)=\frac{\lambda^{2}}{p}$ and 
$\inf_{t \geq 0} h(\lambda, l, t)= 
	\frac{\lambda^{2}}{p} \frac{F\left(\frac{k+p-l}{2}+1, \frac{p}{2}+1, \frac{\lambda^{2}}{2}\right)}
 {F\left(\frac{k+p-l}{2}+1, \frac{p}{2}, \frac{\lambda^{2}}{2}\right)}$;
 \item[ {\bf (b)}] the function $h(\cdot, l, t)$ is, for $t > 0$, $l < k+p$, strictly increasing on 
 $\left[0, \infty\right)$:
 % with $\sup_{\lambda \geq 0} h(\lambda, l, t)= 	\frac{1+t}{t} $;
 \item[ {\bf (c)}]  the function $h(\lambda, \cdot, t)$ is, for $\lambda >0$, $t >0$, strictly increasing 
on $(-\infty, k+p)$; and consequently $h(\lambda, l, t) > h(\lambda, 0, t)$ for $l \in (0,p+k)$.
\end{enumerate}

\end{lemma}
\noindent {\bf Proof.} In part {\bf (a)}, the monotonicity follows from part (a) of Lemma \ref{propertiesratio} setting
$c=0$, and implies $\sup_{t\geq 0}h(\lambda, l, t)= h(\lambda, l, 0)$ and 
$\inf_{t\geq 0}h(\lambda, l, t)=\lim_{t \to \infty} h(\lambda, l, t)$ yielding the results.  
In {\bf (b)}, use the recurrence relation (e.g., Abramowitz and Stegun, 1964) 
$wF(\alpha+1, \beta+1, w)=\beta F(\alpha+1, \beta, w)- \beta F(\alpha, \beta, w)$
% for $a, b, w>0$.
to obtain  
$$h(\lambda, l, t)=\frac{1+t}{t}\left(1-\frac{F(\frac{p+k-l}{2}, \frac{p}{2}, t)}{F(\frac{p+k-l}{2}+1,\frac{p}{2}, t)}\right) \, .$$
The increasing property then follows from part (a) of Lemma \ref{propertiesratio} setting 
$a=\frac{p+k-l}{2} ,b=\frac{p}{2},c= 1$.  Finally, part {\bf (c)} is a consequence of part (b) 
of Lemma \ref{propertiesratio} and expression (\ref{hbul}).  \qed

%For th Bayes estimators with respect to the uniform on the ball, we can conclude that $h^{l}_{U}\left(t\right)\leq %h^{l}_{m}\left(t\right)$ and that $h^{l}_{U}\left(0\right)= m^2/p$, for all $l$.

\begin{remark}
\label{choiceofl}
As shown above in Lemma \ref{propertiesh}, $h^l_{\lambda}(t)$ increases in $l$ so that the amount of shrinkage or expansion of
the boundary uniform estimators is controlled by the choice of the $l$ in the prior (\ref{prior}).  Namely, in comparison to
the benchmark $\delta_{BU,0}$, $\delta_{BU,l}$ will either expand for $k+p > l > 0$, or shrink for $l <0$ with the case $l \to -\infty$ leading to $\delta_{BU,l}$ shrinking to the degenerate estimator $\delta \equiv 0$.   These properties will lead directly to precise risk comparisons in Section 3.  We limit $l$ to be less than $p+k$ in order that the posterior distribution exist for all $(x,s^2)$.  Interestingly for $l \in [p+k, p+k+1)$, the expansion of (\ref{formula}) holds even though the associated distribution is not a probability distribution.
\end{remark}

\subsection{The maximum likelihood estimator}
\noindent
In the same context, we determine the maximum likelihood estimator $\delta_{\hbox{mle}}$ and some of its analytical properties. 

\begin{lemma}\label{mle}
For model (\ref{model}) and parameter space $\Theta(m)$ in (\ref{ps}), the maximum likelihood estimator of $\theta$ is given by
$\delta_{\hbox{mle}}(X,S^2) = h_{\hbox{mle}}(T) X$, with $T=\frac{\left\|X \right\|^{2}}{S^{2}}$ and
\begin{equation}
\label{hmle}
 h_{\hbox{mle}}(t)= \left( \frac{m^{2}}{2(p+k)}\left( \sqrt{1+4\frac{(k+p)}{m^{2}}\frac{1+t}{t}} -1\right)  \wedge 1  \right) 
 \, .
\end{equation}
 
\end{lemma}
\noindent {\bf Proof.}
According to (\ref{model}), the loglikelihood is given by
\begin{eqnarray*}
 \ln L\left(\theta, \sigma\right) & \propto & -\frac{1}{2\sigma^{2}} \left(\left\|x-\theta\right\|^{2}+s^{2}\right)-\frac{p+k}{2} \ln \sigma^{2} \,.
 \end{eqnarray*}
For fixed $\sigma$, the likelihood with respect to $\theta$ is maximized for
$$ \hat{\theta}_{\sigma} =
\left\lbrace \begin{array}{ll} x & \;\hbox{if } \,\;\left\|x\right\| \leq m \sigma  \;  \\ 
   m \sigma \frac{x}{\left\|x\right\|} & \;\,
\hbox{if } \; \left\|x\right\| \geq m \sigma \,.  
\end{array} \right.
$$
Therefore, we have $\sup_{\left(\theta, \sigma\right) \in \Theta\left(m\right)}L\left(\theta, \sigma\right)=\sup_{\sigma > 0}L\left(\hat{\theta}_{\sigma}, \sigma\right)$. Next, we can see that 
\begin{eqnarray*}
 -\ln L\left(\hat{\theta}_{\sigma}, \sigma\right) & \propto & \frac{1}{2\sigma^{2}} \left(\left\|x-\hat{\theta}_{\sigma}\right\|^{2}+s^{2}\right)+\frac{p+k}{2} \ln \sigma^{2}
 \\
 &=& 
\left\lbrace \begin{array}{ll} \frac{s^{2}}{2\sigma^{2}}+ \frac{p+k}{2} \ln \sigma^{2} & \hbox{if } \,\left\|x\right\| \leq m \sigma \; \\ 
   \frac{s^{2}}{2\sigma^{2}}+ \frac{p+k}{2} \ln \sigma^{2}+\frac{1}{2}\left(m^{2}+\frac{\left\|x\right\|^{2}}{\sigma^{2}}-2m\frac{\left\|x\right\|}{\sigma}\right) & \,
\hbox{if }  \left\|x\right\| \geq m \sigma \,.  
\end{array} \right.
\end{eqnarray*}
Now, the minimum of $-\ln L\left(\hat{\theta}_{\sigma}, \sigma\right)$ is attained on $\left(\frac{\left\|x\right\|}{m}, \infty\right)$ whenever $\frac{s^{2}}{p+k}\geq \frac{x' x}{m^{2}}$. Therefore, $\hat{\theta}_{\hbox{mle}}(x,s)=x$ 
whenever $t= \frac{x'x}{s^2} \leq \frac{m^2}{p+k}$.  On the other hand, when $t > \frac{m^2}{p+k}$,
the minimum of $-\ln L\left(\hat{\theta}_{\sigma}, \sigma\right)$ is attained at 
$$\hat{\sigma_{0}}=\frac{m \;\left\|x\right\|}{2(p+k)}\left( \sqrt{1+4\frac{(k+p)}{m^{2}}
\left(\frac{s^{2}}{\left\|x\right\|^{2}}+1\right)} -1\right)\,.$$ 
Finally, since $\hat{\theta}_{\hbox{mle}}(x,s)=m\frac{x}{\left\|x\right\|}\hat{\sigma_{0}}$ for such values of $(x,s)$, 
the result follows.
$\qed$

\begin{remark}
\label{mleshrinks}
Observe that 
$\|\delta_{\hbox{mle}}(x, s^2) \| \leq \left\|x\right\|$ with strict inequality for $\frac{x'x}{s^2} > \frac{m^2}{p+k}$, so that the maximum likelihood
shrinks the unbiased estimator towards the origin.  Moreover, the amplitude of this shrinkage as measured by 
$h_{\hbox{mle}}(t)$ is seen to be increasing with respect to $t$ (i.e., $h_{\hbox{mle}}(t)$ decreases in 
$t$), with $\lim_{t \rightarrow \infty} h_{mle}(t)=\frac{\sqrt{1+4\gamma}-1}{2\gamma} \,,$
$\gamma=\frac{p+k}{m^{2}}$.  
\end{remark}

Now, for various risk differences involving the estimators $\delta_{\hbox{mle}}$, $\delta_{UB}$, and the Bayes estimators
$\delta_{BU,\,l}$, the following comparisons will be pivotal.
 
\begin{lemma}\label{tough}   
\begin{enumerate}
\item[{\bf (a)}]  We have $h(m,0,t) \leq h_{\hbox{mle}}(t)$ for
$t\geq \frac{m^{2}}{p+k}$, whenever $p \geq 2, k \geq 2$, and whenever $k=1$ and $m \leq \sqrt{p}$.
\item[{\bf (b)}]  For $p \geq 2$ and for $p=1, k=1$, the condition $m \leq \sqrt{p}$ is necessary and sufficient so that 
$h(m,0,t) \leq h_{\hbox{mle}}(t)$ for all 
$t \geq 0$ and, in such cases, equality is attained if and only if $t=0$ and $m = \sqrt{p}$.
\item[{\bf (c)}]  Whenever $m \leq \sqrt{p}$ and $l < p+k$, we have $h(m,l,t) \leq 1$ for all 
$t \geq 0$ with equality if and only if $t=0$ and $m = \sqrt{p}$.
\end{enumerate}
\end{lemma}
\noindent {\bf Proof.}  Part  {\bf (c)} is a direct consequence of part (a) of Lemma \ref{propertiesh}.
Since $h_{\hbox{mle}}(\cdot)$ is constant on $[0, \frac{m^{2}}{p+k}]$, and 
$h(m,0,\cdot)$ is decreasing and bounded above by $\frac{m^2}{p}$ on $[0, \frac{m^{2}}{p+k}]$ by virtue of
Lemma \ref{propertiesh}, part {\bf (a)} implies  part {\bf (b)}.
For proving part {\bf (a)}, by making use of representations (\ref{hmle}) and (\ref{hbul}), it suffices to verify that
\begin{eqnarray} 
\nonumber  \frac{F\left(\frac{k+p}{2}+1, \frac{p}{2}+1, z \right)}
 {F\left(\frac{k+p}{2}+1, \frac{p}{2}, z \right)}
 & \leq &
 \frac{p}{2(p+k)}\left( \sqrt{1+\frac{2(k+p)}{ z}} -1\right)\,, \\
\label{R1} R(z)=  \frac{F\left(\frac{k+p}{2}+1, \frac{p}{2}, z \right)}
 {F\left(\frac{k+p}{2}+1, \frac{p}{2}+1, z \right)}
 & \geq &  
 \frac{z}{p}\left( \sqrt{1+\frac{2(k+p)}{ z}} +1\right)\,,
\end{eqnarray}
for all $p \geq 2, k \geq 2, z>0$ or for $k=1$, $m \leq \sqrt{p}$ and $z=\frac{m^{2} t}{2(1+t)} \in [\frac{m^{4}}{2(m^{2}+p+k)}, \frac{m^2}{2})$. 
With the decomposition 
\begin{eqnarray*}
F\left(\frac{k+p}{2}+1, \frac{p}{2}, z \right) &=& \sum_{j \geq 0} \frac{(\frac{k+p}{2}+1)_j}{(\frac{p}{2}+1)_j}
\, \frac{(\frac{p}{2}+1)_j}{(\frac{p}{2})_j} \, \frac{z^j}{j!} = \sum_{j \geq 0} \frac{(\frac{k+p}{2}+1)_j}{(\frac{p}{2}+1)_j}
\, (1 + \frac{2j}{p}) \, \frac{z^j}{j!} \\
\, &=& F\left(\frac{k+p}{2}+1, \frac{p}{2} + 1, z \right) + \frac{2z}{p} \sum_{j \geq 0} 
\frac{(\frac{k+p}{2}+1)_{j+1}}{(\frac{p}{2}+1)_{j+1}} \, \frac{z^j}{j!}\,, 
%\\
%\, &=& F\left(\frac{k+p}{2}+1, \frac{p}{2} + 1, z \right) + \frac{2z}{p} (1 + \beta(z))
\end{eqnarray*}
one obtains the representation 
\begin{equation}
\label{R2}
 R(z) = 1 + \frac{2z}{p} \beta(z) = 1+ \frac{2z}{p} (1 + k E_z[\frac{1}{p+2J+2}]),
\end{equation} 
with $E_z$ representing expected value with respect to a discrete random variable $J$ with probability mass function
proportional to $  \frac{(\frac{k+p}{2}+1)_{j}}{(\frac{p}{2}+1)_{j}} \frac{z^j}{j!} \, 1_{\{0,1, \ldots\}}$.
With a similar expansion, we also have 
\begin{equation}
\label{zbeta}
 z (\beta(z)-1) = k E_z[\frac{J}{k+p+2J}]\,.  
\end{equation} 
   Using (\ref{R2}) and (\ref{zbeta}), (\ref{R1})  becomes equivalent to
\begin{eqnarray}
\nonumber  & 1+ \frac{2z}{p} \beta(z)  \geq  \frac{z}{p} \, (\sqrt{1+\frac{2(k+p)}{ z}} +1 ) \, \\
\label{A>}
\, \Longleftrightarrow & A(z)= \frac{1}{2k} \left\lbrace 4 \beta(z) (\beta(z) - 1) z^2 + 4 p (\beta(z) - 1) z + p^2
\right\rbrace \geq z \\
\label{sc}
\, \Longleftrightarrow & \frac{p^2}{2kz} + T(z) \geq 1\,,
\end{eqnarray}
where 
\begin{equation}
T(z)= \frac{2}{k} \{\beta(z) (\beta(z) - 1) z + p (\beta(z) - 1) \}\,.
\end{equation}
Now, by making use of (\ref{R2}) and (\ref{zbeta}), expand $T(\cdot)$ in terms of two independent copies $J_1$ and $J_2$ of $2J$ as
\begin{eqnarray*}
\nonumber   T(z)&=&\left( 1+ k E_z(\frac{1}{p+2+J_1})  \right) \left(E_z(\frac{J_2}{k+p+J_2}) \right) + E_z(\frac{2p}{p+2+J_1}) \\
\,  &=&  E_z\left[\frac{2p(k+p)  + J_1J_2 + (p+2) J_2 + (k+2p) J_2}{(p+2)(k+p)  + J_1J_2 + (p+2) J_2 + (k+p) J_1}  \right].
\end{eqnarray*}
In the above, we have expressed the product of expectations as the expectation of a product given that the $J_i$'s are independent.
Given (\ref{sc}), it will suffice to show that $T(z)$ is lower bounded by 1 for all $z >0$.  Since $p \geq 2$ by assumption, we obtain with 
an expansion where we denote $p$ the joint probability function of $(J_1,J_2)$ and $d(j_1,j_2)=(p+2)(k+p)  + j_1j_2 + (p+2) j_2 + (k+p) j_1$,
\begin{eqnarray*}
\nonumber   T(z)& \geq&  E_z\left[\frac{(p+2)(k+p)  + J_1J_2 + (p+2) J_2 + (k+p) J_2}{(p+2)(k+p)  + J_1J_2 + (p+2) J_2 + (k+p) J_1}  \right] \\
\, &= & 1 + (k+p) E_z\left[\frac{J_2-J_1}{(p+2)(k+p)  + J_1J_2 + (p+2) J_2 + (k+p) J_1}  \right] \\
\, &= & 1 + (k+p) (\sum_{j_2 < j_1} + \sum_{j_2 > j_1}) \left(\frac{(j_2 - j_1) \, p(j_1,j_2)}{d(j_1,j_2)}\right) \\
\, &= &  1 + (k+p)  (\sum_{j_2 < j_1} (j_2-j_1) \left(\frac{p(j_1,j_2)}{d(j_1,j_2)} -  \frac{p(j_2,j_1)}{d(j_2,j_1)} \right) \\
\, &=& 1 + (k+p) (k-2) \sum_{j_2 < j_1} p(j_1,j_2) \frac{(j_2-j_1)^2}{d(j_2,j_1) \,d(j_1,j_2)} \\
\, &\geq& 1\,,
\end{eqnarray*}
since $J_1$ and $J_2$ are independent and for $k \geq 2$.  This establishes the result for all $(m,p,k)$ such that $p \geq 2$ and $k \geq 2$.  Finally, for $k=1$ and $m \leq \sqrt{p}$, it is clear that (\ref{A>}) is verified for $z \leq m^2/2$.  This completes the proof. \qed
\footnote{We point out that we have shown inequality (\ref{R1}) for all
$p \geq 2, k \geq 2, z \geq 0$. We could not find a similar or as sharp as inequality in the literature.  Such an inequality, and the technique used to establish it, may well be of independent interest as various bounds for ratios of special functions like the hypergeometric here have been previously the subject of study (e.g., Joshi and Bissu, 1996; Kokologiannaki, 2012.)
}

\section{Dominance results }

Following Lemma \ref{risklemma}'s representation of risks, and the analytical properties of $h(m,l,\cdot)$ and $h_{\hbox{mle}}(\cdot)$
worked out in Lemmas \ref{propertiesh} and \ref{tough}, we obtain various risk comparisons and dominance results which
will be consequences of the following.

\begin{theorem}
\label{riskcomparison}
Consider equivariant estimators $\delta_{h}$ as in (\ref{deltah}) for 
estimating $\theta$ for model (\ref{model}), loss (\ref{loss}) and the parameter space $\Theta(m)$. 
\begin{enumerate}
\item[ \bf{(a)}] A sufficient condition for $\delta_{h_2}$ to dominate $\delta_{h_1}$ ($\delta_{h_2} \neq \delta_{h_1}$) is $h_1(t) \geq h_2(t)$ 
and $\frac{h_{2}(t) +  h_{1}(t)}{2} \geq h(m,0,t)$ for all $t>0$;
\item[ \bf{(b)}]  Let $S=\{t: h_1(t) > h(m,0,t)\}$.  Let $\underline{S}$ be a subset (not strict) of $S$ such that  
$\nu(\underline{S}) > 0$, where  $\nu$
is the Lebesgue measure on $\mathbb{R}$, then
$\delta_{h_1}$ is inadmissible, and dominated by $\delta_{h}$ with $h(t) = h_1(t) \wedge h(m,0,t)$, as well as
by any $\delta_{h}$, $h(t)=h_1(t) 1_{\underline{S}^c}(t) + h^*(t) 1_{\underline{S}}(t)$ with $ 2h(m,0,t) - h_1(t) \leq h^*(t) \leq h_1(t)$
for all $t \in \underline{S}$.
\item[ \bf{(c)}]  If $h_1(t) \geq h(m,0,t)$ for all $t \geq 0$, $h_1(\cdot) \neq h(m,0,\cdot)$, then the boundary uniform Bayes estimator $\delta_{BU, \, 0}$ given in (\ref{hbul}) dominates $\delta_{h_1}$. 
\end{enumerate}
\end{theorem}
{\bf Proof.}  Parts {\bf (b)} and {\bf (c)} are consequences of part {\bf (a)} by setting $h_2 \equiv h$ (as given) and
$h_2(\cdot) \equiv h(m,0,\cdot)$ respectively.  From (\ref{riskdecomposition}), we obtain
\begin{eqnarray} \label{7}
	R(\lambda, \delta_{h_{1}}) - R(\lambda, \delta_{h_{2}}) &=& 
	E_{\lambda} \left[ a(T) \{\left(h_{1}(T)-h_{2}(T)\right) \times \left(h_{1}(T)+ h_{2}(T)- 2 h(\lambda,0,T\right)\} \right] \\ &=&	 E_{\lambda}\left[a(T) \; \Delta_{\lambda}\left(T\right)\right] \; (\hbox{say})\,,
\end{eqnarray}
where $a(\cdot)(> 0)$ is given in Lemma \ref{risklemma}.  Therefore, the conditions on $h_1$ and $h_2$, along with part (b) of Lemma \ref{propertiesh} which implies $h(\lambda,0,\cdot) < h(m,0,\cdot)$ for $\lambda < m$, force $\Delta_{\lambda}(T) \geq 0$ with probability one for all $\lambda \in [0,m]$.
Finally,  dominance indeed occurs with $\Delta_{\lambda}(t) > 0$ for all $\lambda \in [0,m)$ and $t \in \{u: h_1(u) > h_2(u)\}$.   $ \; \; \; \; \; \; \; \;\; \; \; \; \; \qed$

\begin{corollary}\label{dominanceresults}
\begin{enumerate}
\item[ {\bf (a)}] For $\left(m, p,k\right)$ such that $p \geq 2, k \geq 2$ or $k=1, m\leq \sqrt{p}$, 
$\delta_{\hbox{mle}}$ dominates the unbiased estimator $X $. 
\item[ {\bf (b)}]  For $m\leq \sqrt{p}$, $\delta_{BU,\,0}$ dominates $\delta_{\hbox{mle}}$ whenever $p \geq 2$ or
$(p,k)=(1,1)$. 
\item[ {\bf (c)}]  Whenever $m\leq \sqrt{p}$ and $l \in [0,k+p)$, $\delta_{BU,\,l}$ dominates $X$. 
\item[ {\bf (d)}]  For all $(m,p)$ and $l \in (0,k+p)$, $\delta_{BU,\,l}$ is inadmissible and dominated by
$\delta_{BU,\,0}$.
\item[ {\bf (e)}]  For all $(m,p)$ with $m>\sqrt{p}$, such that $p \geq 2, k \geq 2$, the estimator $\delta_{\hbox{mle}}$ is inadmissible and dominated by the truncation $\delta_h$, with $h(t)=h_{\hbox{mle}}(t) \wedge h(m,0,t)$. \footnote{In both {\bf (d)} and {\bf (e)}, further dominating procedures can be derived using part {(b)} of Theorem \ref{riskcomparison}}
\end{enumerate}
\end{corollary}
\noindent{\bf Proof.} {\bf (a)}  We apply part (b) of Theorem \ref{riskcomparison} with $h_1 \equiv 1$, 
$\underline{S}=(\frac{m^2}{p+k}, \infty)$, and $h^* \equiv h_{\hbox{mle}}$.  Recall that $h_{\hbox{mle}}(t) < 1$ 
for $t > \frac{m^2}{p+k}$ as pointed out 
in Remark \ref{mleshrinks}.  Under the given conditions on $\left(m, p,k\right)$, part (a) of Lemma \ref{tough} tells us that $S=\{t: 1 > h(m,0,t)\}$ contains $\underline{S}$ and that $h(m,0,t) \leq h_{\hbox{mle}}(t) < h_1(t)$ for all $t \in \underline{S}$, from which the result follows. \\
{\bf (b)}  Under the assumptions on $\left(m, p,k\right)$, part (b) of Lemma \ref{tough} tells us that $h_{\hbox{mle}}(\cdot) \geq h(m,0,\cdot)$ and the result is hence an immediate consequence of part (c) of
Theorem \ref{riskcomparison}. \\
{\bf (c)} Parts (c) of Lemma \ref{propertiesh} and Lemma \ref{tough} imply that $h(m,0,t) < h(m,l,t) < 1$ for all $t>0$ in cases where $m\leq \sqrt{p}$ and $l \in [0, p+k)$.  The result then follows as an application of part (b) of Theorem \ref{riskcomparison} with $h_1 \equiv 1$, $h^*(\cdot)=h(m,l,\cdot)$, and $\underline{S}=S=(0,\infty)$. \\
{\bf (d)}  Given part (c) of Lemma \ref{propertiesh}, the result follows as an application of part (c) of Theorem \ref{riskcomparison}. \\
{\bf (e)}  Follows from part (b) of Theorem \ref{riskcomparison} with $h_1(\cdot)= h_{\hbox{mle}}(\cdot)$, $h^*(\cdot)=h(m,0,\cdot)$ and $\underline{S}=S=(\frac{m^2}{p+k},\infty)$, and by making use of part (a) of Lemma \ref{tough}. 
\qed

\begin{remark}
The dominance results in part (c) of Theorem \ref{dominanceresults} were obtained by Kubokawa (2005) in the univariate case,
and constitute an extension to the multivariate case $p>1$ for improving on $\delta_{UB}$.  However, with part (a) demonstrating the superiority of $\delta_{\hbox{mle}}$ on  $\delta_{UB}$, the more important finding is the dominance
of $\delta_{BU,\,0}$ over $\delta_{\hbox{mle}}$ for $m \leq \sqrt{p}$ and $p \geq 2$.\footnote{For completeness, we have also summarized the situation in the very special case $p=k=1$.} For cases where $m > \sqrt{p}$, $p,k \geq 2$, $\delta_{\hbox{mle}}$ is still inadmissible and explicit improvements are available from part (e).
  We believe and conjecture that these results extend to the univariate case, but we have been unable to establish the inequality in part (a) of Lemma \ref{tough}, which is a critical element of the analysis.   
\end{remark}

\begin{example}
\label{projection}
Part (b) of Theorem \ref{riskcomparison} indicates clearly that equivariant estimators $\delta_h$ with $h$ taking values that are too large on a subset of $\mathbb{R}^+$ are inefficient and can be improved upon by projecting towards the benchmark 
$h(m,0,\cdot)$.  Otherwise said, the function $h(m,0,\cdot)$ provides an upper envelope for a complete class of estimators.  Key applications are explicited in Corollary \ref{dominanceresults}, leading to improvements on the unbiased estimator, the maximum likelihood estimator, and the generalized Bayes estimators $\delta_{BU,l}$ with $0<l < p+k$.  Moreover, the degree of expansion translates to even greater inefficiency as is the case for the unbiased estimator.  As an example, taking $(p,k,m)=(5,10,2)$, Figure 1 shows the multipliers and risk functions of $\delta_{UB}$, 
$\delta_{\hbox{mle}}$ and $\delta_{BU,0}$ as a function of $\lambda= \frac{\|\theta \|}{\sigma} \in [0,2]$.  Here $m \leq \sqrt{p}$, so that the ordering of the multipliers and risks is clearly dictated, namely by Lemma \ref{tough} and Corollary \ref{dominanceresults}.  We see that $h_{UB}=1$ is much too large, with very poor risk performance.  The estimator $\delta_{\hbox{mle}}$ fares better but the gains provided by the Bayes estimator $\delta_{BU,0}$ are nevertheless important; in relative terms ranging from a maximum of over 50\% at $\lambda=0$ to a minimum of around 12\% at the boundary $\lambda=2$.
\end{example}

%   \begin{figure}[ht]
%    \centering
%  \includegraphics[width= 10cm]{figure1.pdf}
%  \caption{Graphs of the ratio of related improvement in risks for $d=2$ with $a=0.71, a_{\hbox{tlx}}=0.8, a=0.85$ and  %$a=0.90$.}\label{figure1}
%   \end{figure}\renewcommand{\theequation}{\thesection.\arabic{equation}}
%   \begin{figure}[ht]
%    \centering
%  \includegraphics[width= 10cm]{figure2.pdf}
%  \caption{Graphs of the ratio of related improvement in risks for $d=2$ with $a=0.71, a_{\hbox{tlx}}=0.8, a=0.85$ and  %$a=0.90$.}\label{figure1}
%   \end{figure}\renewcommand{\theequation}{\thesection.\arabic{equation}}
   
 \begin{figure}[ht]
 \centering
 \includegraphics[width=7.5cm]{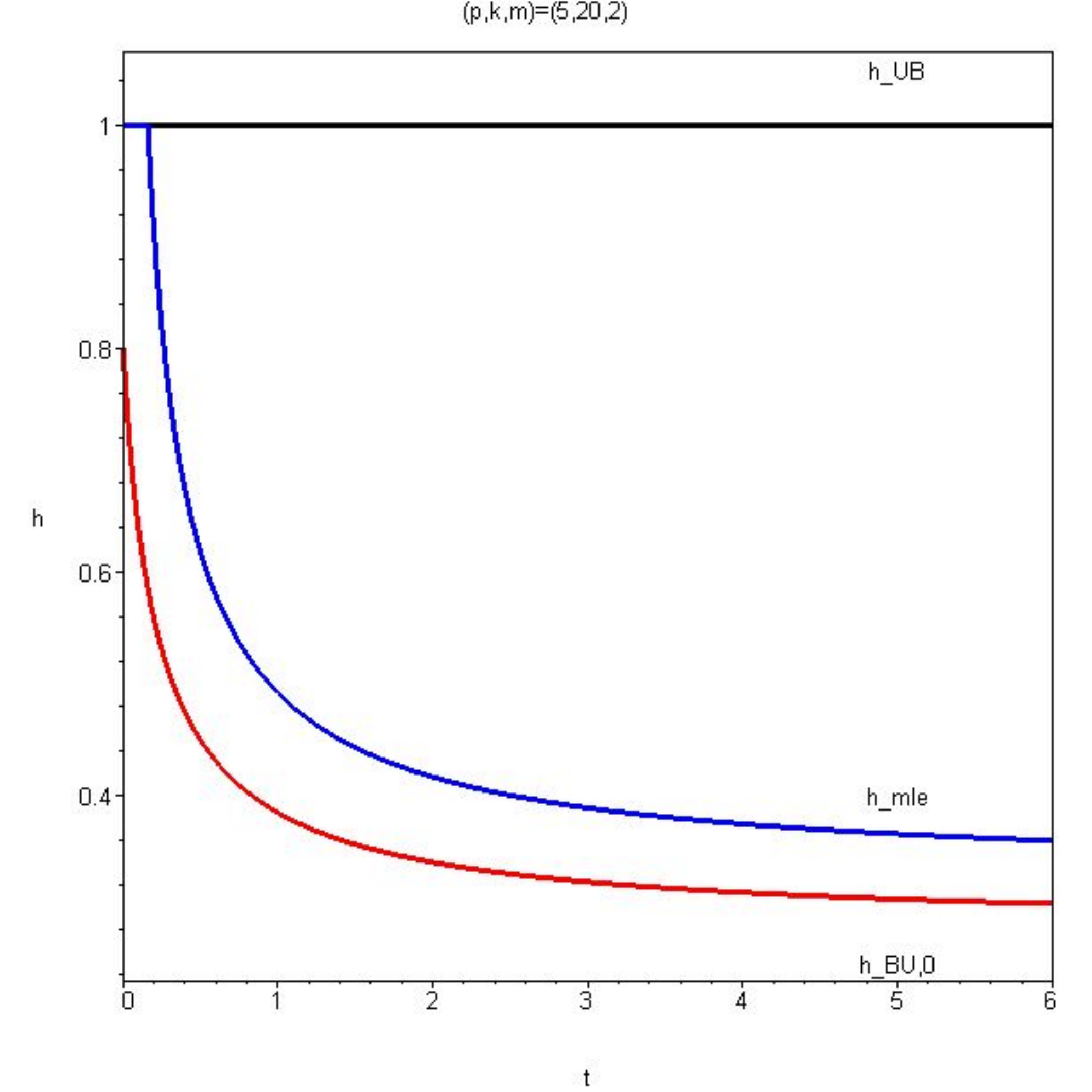}
 \hspace{1cm}
 \includegraphics[width=7.5cm]{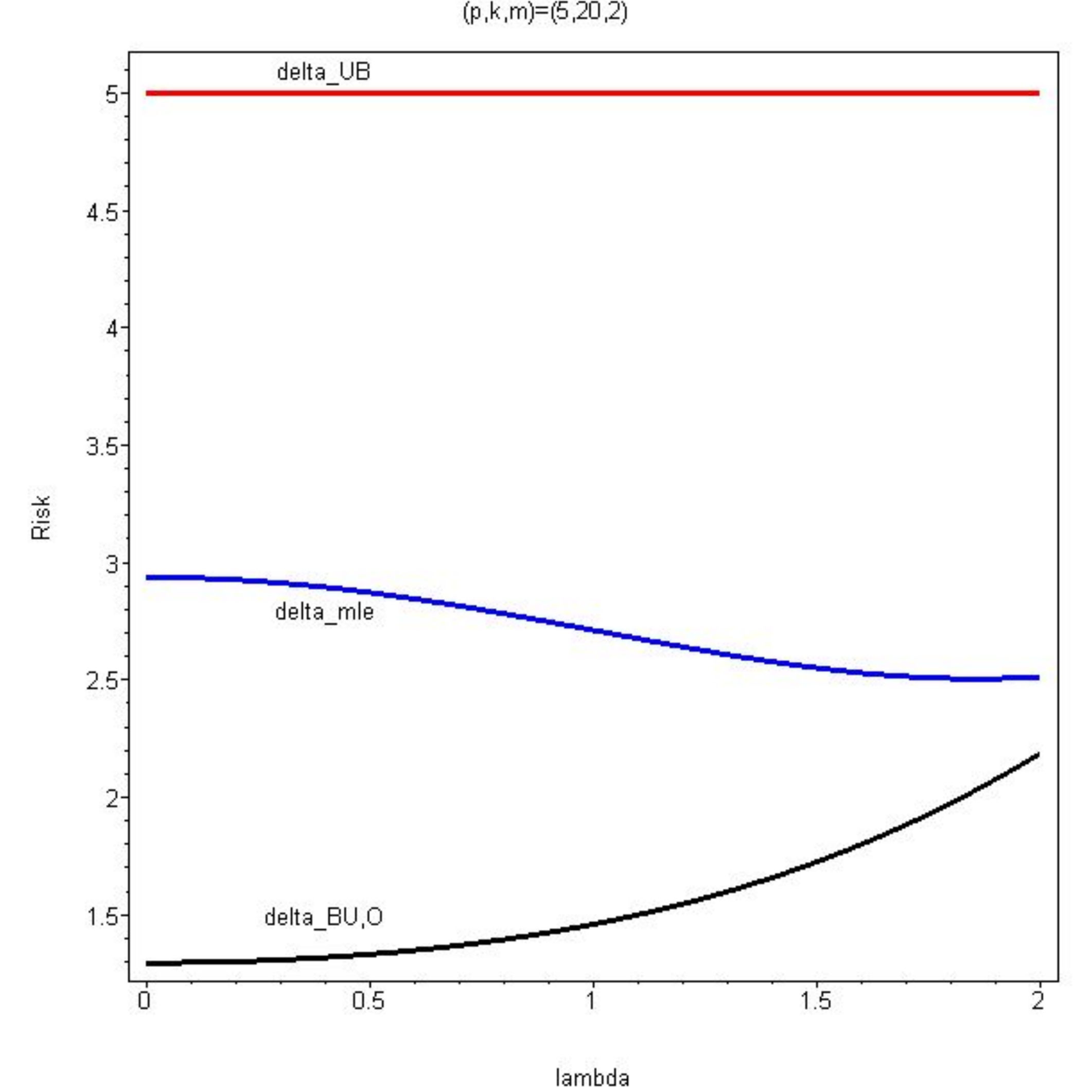}
  \caption{Multipliers and risks of $\delta_{UB}$, $\delta_{\hbox{mle}}$ and $\delta_{BU,0}$ for $(p,k,m)=(5,20, 2)$. }\label{graphe_001}
 \end{figure}

\begin{example}
\label{example}
Further risk function comparisons are presented in Figure 2 for other combinations of $(p,k,m)$.  In both cases, we do not have
$m \leq \sqrt{p}$ and the dominance findings for the boundary Bayes estimator $\delta_{BU,0}$ do not apply, and there is also no guarantee that it performs satisfactorily from a frequentist risk point of view, even in comparison to the unbiased estimator.
In one of the cases (with $p=5$ and $m=3$), the numerical results indicate that $\delta_{BU,0}$ performs quite well in comparison to $\delta_{\hbox{mle}}$ with dominance and relative gains between 5\% and 20\%.  But notice how the gains are less impressive than in Figure 1 where $m$ is smaller and the dimension $p=5$ is the same.  In the other case, where the ratio of $m$ relative to $p$ is much larger, not only do the findings not apply, notwithstanding the dominance result applicable to the truncation given in part (e) of Corollary \ref{dominanceresults}, but the risk performance of $\delta_{BU,0}$ is arguably quite poor.  More research is thus required on alternative Bayes estimators, especially when $m > \sqrt{p}$.  The deficiency of $\delta_{BU,0}$ lies in the fact that it expands too much.  Here, as an example, we have $m^2/p=(3^2)/3=3$ so that 
$\delta_{BU,0}(x,s^2) \approx 3 x$ when $x$ is in a neighbourhood of $0$, which leads to poor estimates when $\|\theta\|/\sigma$ is small.  Other priors with $l \leq 0$ will shrink $\delta_{BU,0}$ towards the origin.  One such class of choices, studied in the known variance case by Fourdrinier and Marchand (2010), are obtained by taking $\pi_{\sigma}$ in (\ref{generalprior}) to be 
a uniform density on the sphere of radius $\alpha \sigma$, for all $\sigma >0$ with $\alpha \in [0,m)$.   Of course, for such priors, Theorem \ref{hbul} provides an expression for the Bayes estimator by replacing $m$ by $\alpha$.   Here,  we tried one such choice with $\alpha=2.5$ and the numerical evaluation indicates quite satisfactory performance with a significant  improvement on $\delta_{BU,0}$ near the centre and on a large part of the parameter space, with a slightly worse performance near or on the boundary.   
\end{example}

 \begin{figure}[ht]
 \centering
 \includegraphics[width=7.5cm]{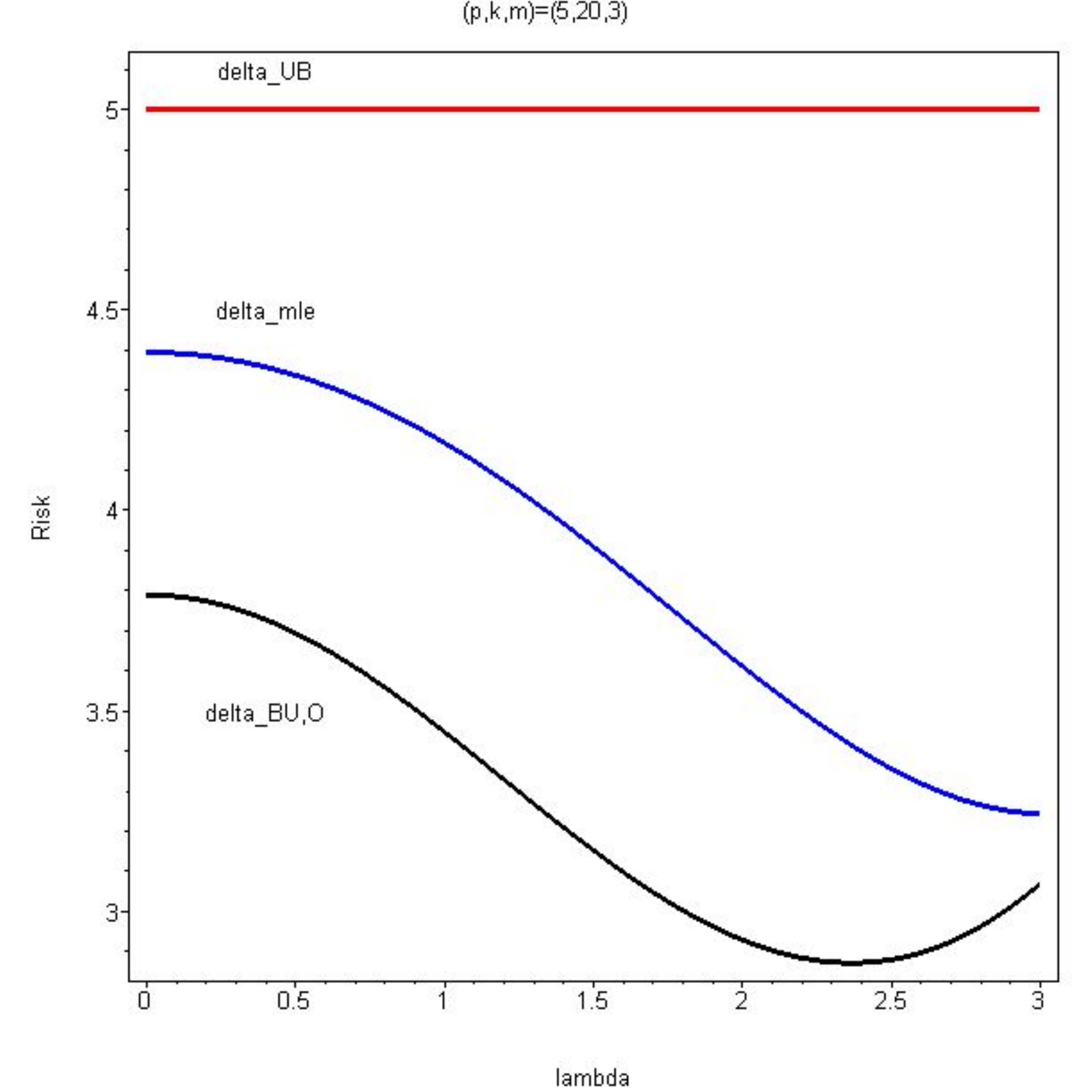}
 \hspace{1cm}
 \includegraphics[width=7.5cm]{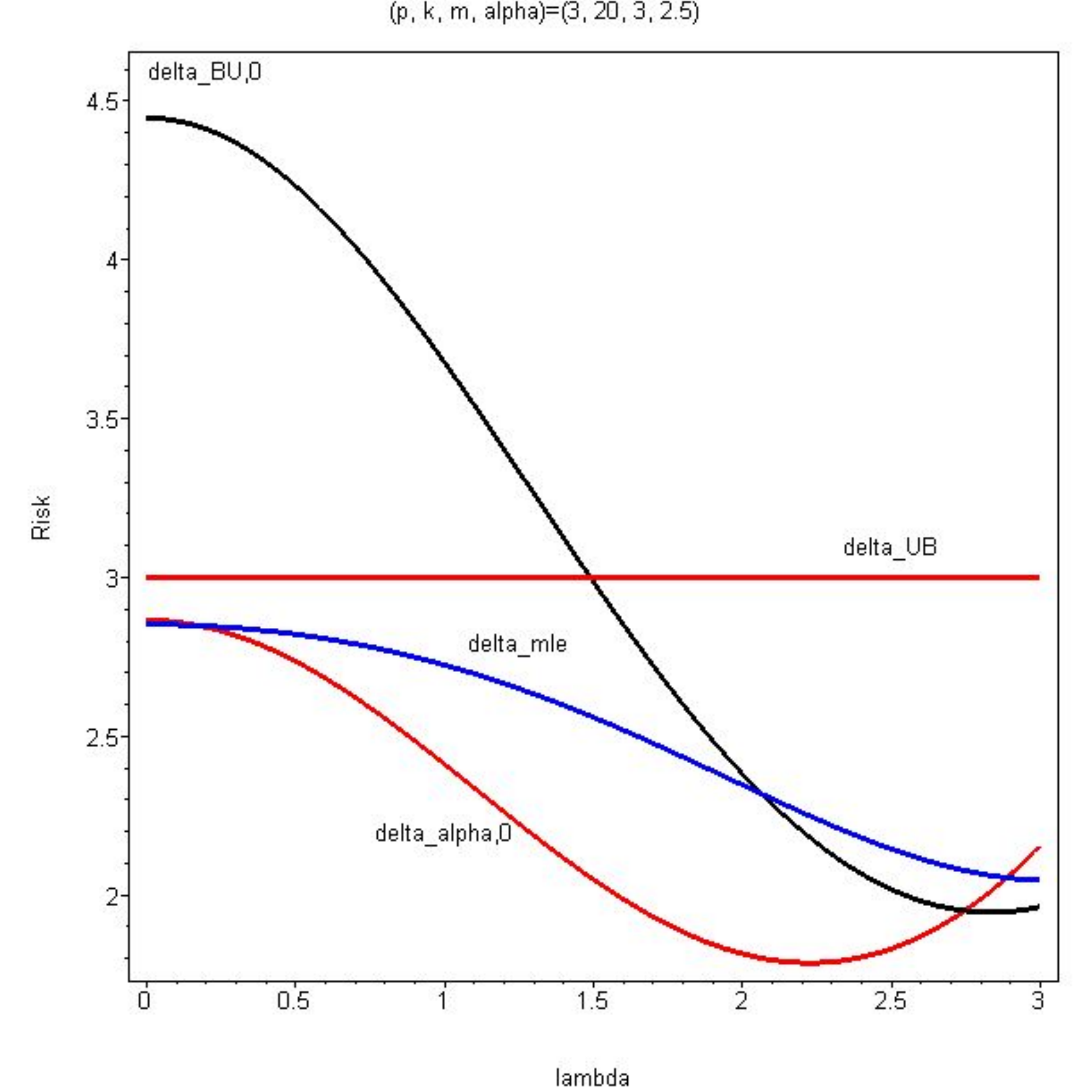}
  \caption{Risks of various estimators : unbiased ($\delta_{UB}$), maximum likelihood ($\delta_{\hbox{mle}}$), boundary uniform ($\delta_{BU,0}$) and Bayes with $\theta|\sigma$ uniform on sphere of radius $\alpha \sigma$ ($\delta_{\alpha,0}$), for $(p,k,m)=(5,20, 3)$ and $(p,k,m,\alpha)=(3,20, 3,2.5)$}\label{graphe_001}
 \end{figure}

We conclude with a general representation for Bayesian estimators and a universal dominance result applicable to very small parameter spaces (where the average squared signal to noise ratio is less than or equal to $1/2$) which focusses quite clearly on the inadequacy of the unbiased estimator and lies in continuity with its inefficiency among affine linear estimators described in the introduction and the above risk comparisons.

\begin{lemma}
\label{ss}
For priors as in (\ref{generalprior}) with spherically symmetric densities $\pi_{\sigma}$ for all $\sigma$, Bayes estimators are equivariant, of the form $\delta_{\pi}(x,s^2)= h_{\pi}(t) \, x$, with {\bf (i)} $0 \leq h_{\pi}(t) \leq h(m,l,t)$ for all $t >0$, and with {\bf (ii)} $ h_{\pi}(\cdot) \leq h(m,0,\cdot)$ for all $l < 0$.
\end{lemma}
{\bf Proof.}  Part {\bf (ii)} is immediate from part {\bf (i)} and part (b) of Lemma \ref{propertiesh}.  For {\bf (i)},  first observe that Bayes estimators in (\ref{formula}) are also posterior expectations $E_{\pi^*}(\theta|x,s^2)$ for prior measures $\pi^*(\theta, \sigma)= \pi_{\sigma}(\theta) \sigma^{l-4}$ (as expressed in equation \ref{ar}).  Now, proceeding as in Marchand and Perron (2001, Theorem 4), with spherically symmetric choices, the prior densities $\pi_{\sigma}$ in (\ref{generalprior}) admit the representation $\theta|\sigma =^d RU$, where $R$ is supported on $[0,m]$, $U$ is uniformly distributed on the sphere $S_{\sigma}$, and $R$ and $U$ are independent (conditional on $\sigma$).  Hence, by making use of the developments in Theorem \ref{bul}, we have
$\delta_{\pi}(x,s^2)=E_{\pi^*}(\theta|x,s^2)=E[R E(U|R,x,s^2)]= E[h(R,l,t)|x,s^2] \, x$, where $h(r,l,t)$ is given in (\ref{hbul})
and $R|x,s^2$ is the posterior distribution of $\|\theta \|$.  
Since this posterior distribution is supported on $[0,m]$, the result follows since $h(\cdot, l, t)$ is nonnegative and increasing for all $(l,t)$ by virtue of part (a) of Lemma \ref{propertiesh}. \qed 

\begin{theorem}
\label{universal}  
Whenever $m \leq \sqrt{\frac{p}{2}}$, all equivariant estimators $\delta_h(x,s^2)= h(t) x$ with $0 \leq h(t) \leq h(m,0,t)$ for all $t>0$ dominate the unbiased estimator $X$.  These include all Bayes estimators with respect to a prior $\pi$ as in (\ref{generalprior}) with spherically symmetric densities $\pi_{\sigma}$ for all $\sigma$ and $l \leq 0$.
\end{theorem}

{\bf Proof.}  Applying part (a) of Theorem \ref{riskcomparison} with $h_1 \equiv 1$ and $h_2 \equiv h$ with the given assumptions on $h$, $m$ and $p$, we have $$ \frac{h(t) + 1}{2} \geq \frac{1}{2} \geq \frac{m^2}{p} \geq h(m,0,t)$$
for all $t>0$, with the rightmost inequality a consequence of part (a) of Lemma \ref{propertiesh}.
The proof is complete by observing that the inclusion of the Bayesian estimators among the dominating estimators is a consequence of Lemma \ref{ss}.  \qed

\begin{remark} 
The same result holds if the variance $\sigma$ is known.  This corresponds to the setting studied by Marchand and Perron (2001) and the result may be established along the same lines and using similar developments given in their paper.  They actually establish a universal dominance result (also see Fourdrinier and Marchand, 2010) but it applies for the more challenging problem of dominating the maximum likelihood estimator.  All in all, the result here and the other dominance results involving the unbiased estimator may not be that surprising given that even the trivial estimator $\delta \equiv 0$ dominates $X$ whenever $m \leq \sqrt{p}$.  But, we have still provided a method of proof applicable to a very large class of Bayesian estimators.  And, our findings elsewhere relate as well to $\delta_{\hbox{mle}}$.
\end{remark}

\section{Concluding Remarks }

For estimating a multivariate normal mean with an upper bounded signal to noise ratio $\frac{\|\theta\|}{\sigma}$,
we have provided dominance results which can viewed as both multivariate extensions of results
obtained by Kubokawa (2005), as well unknown variance extensions of results obtained by Marchand and Perron (2001).  
In opposition to similar extensions for Stein estimation, the presence of the unknown scale here leads to challenges in describing Bayes estimators and some of their analytical properties which ultimately relate to frequentist risk performance. 

We have focussed mostly on boundary Bayes estimators and the benchmark estimators that are improved upon are those obtained from the principles of unbiasedness and maximum likelihood.  As illustrated theoretically and numerically, and analogously to the known $\sigma$ case, the relative merits of the boundary Bayes procedure seem to fairly well correlate with the ratio of the radius $m$ relative to $\sqrt{p}$ (see as well Marchand and Perron, 2001; Fourdrinier and Marchand, 2010; Kortbi and Marchand, 2012).  

More research is required to assess the performance of other Bayes estimators and namely to propose more attractive choices when $m$ is larger relative to $\sqrt{p}$.  Such alternatives include fully uniform Bayes estimators.  In this regard, an interesting question is whether the estimator $\delta_{U,0}$ improves for all $(m,p)$ on $X$ with an affirmative answer representing an unknown variance extension of Hartigan's (2004) result in the particular case of balls.   

Several other related problems or issues are also of interest.  As, an example, our findings do not address directly the  issues of minimaxity and admissibility of Bayesian estimators.  But it seems plausible and we conjecture that $\delta_{BU,0}$
is minimax for small enough $m$ as suggested by its risk in Figure 1 with the maximum risk attained on the boundary and since
$\delta_{BU,0}$ is quite likely a candidate to be an extended Bayes procedure.  

Finally, we point out that the results obtained here are applicable to two-sample problems with additional information as described by Marchand and Strawderman (2004) and Marchand, Jafari Jozani and Tripathi, 2012.  These involve independently distributed observables $X_1 \sim N_p(\theta_1, \sigma'^2), X_2 \sim N_p(\theta_2, \sigma'^2), S'^2 \sim \sigma'^2 \chi^{2}_{k}$, and the objective of estimating $\theta_1$ with the additional information that $\frac{\|\theta_1 - \theta_2 \|}{\sigma'} \leq m'$.
To achieve this, one ''rotates'' $X_1$ and $X_2$ to the independent coordinates $X=(X_1-X_2)/2$ and $W=(X_1+X_2)/2$
and considers estimators of $\theta_1$ of the form $\delta_{\psi}(X, W, S'^2)=W + \psi(X,S'^2)$ showing that 
$\delta_{\psi_1}$ dominates $\delta_{\psi_2}$ for estimating $\theta_1$ under loss $\frac{\|\delta-\theta_1\|^2}{\sigma^2}$ with the additional information that $\frac{\|\theta_1 - \theta_2 \|}{\sigma} \leq m'$ {\bf if and only if} $\psi_1(X, S^2)$ dominates $\psi_2(X, S^2)$ for estimating $\theta=E(X)=\frac{\theta_1-\theta_2}{2}$ under loss $\frac{\|\psi -\theta\|^2}{\sigma^2}$ and constraint $\frac{\|\theta\|}{\sigma} \leq m$, for model (\ref{model}) with $X \sim N_p(\theta, 
\sigma^2=\frac{\sigma'^2}{2})$, $S^2=\frac{S'^2}{2}$ and $m=\frac{m'}{\sqrt{2}}$.

\subsection*{Acknowledgments}
The research work of \'Eric Marchand is partially supported by NSERC of Canada.
During Othmane Kortbi's Ph.D. studies at the Universit\'e de Sherbrooke, he benefited
from financial support from several sources but he wishes to thank especially the
ISM (Institut de sciences math\'ematiques) and the CRM (Centre de recherches math\'ematiques).
Finally, the authors are grateful to Bill Strawderman and Dominique Fourdrinier for useful discussions
and encouraging us to pursue work on this problem.

\subsection*{References}
\small
\renewcommand{\baselinestretch}{1.2}
\begin{description}

\item  Abramowitz, M. and Stegun, I. (1966).  {\it Handbook of mathematical functions
with formulas, graphs, and mathematical tables}.  Dover, New York.

%\item  Berger, J.O. (1985).  {\em Statistical Decision Theory and
%Bayesian Analysis}, Springer Texts in Statistics, Springer-Verlag,
%second edition, New York.

\item  Eaton, M. L. (1989).  {\it Group invariance applications in statistics.} Regional 
Conference Series in Probability and Statistics, vol. 1, Institute of
Mathematical Statistics, Hayward, California. 

\item Fourdrinier, D. \& Marchand, \'E. (2010).
On {B}ayes estimators with uniform
priors on spheres and their comparative performance with maximum likelihood
estimators for estimating bounded multivariate normal means.
{\it Journal of Multivariate Analysis},
{\bf 101}, 1390-1399.  

\item Hartigan, J. (2004).  Uniform priors on convex sets improve
risk. {\it Statistics \& Probability Letters}, {\bf 67}, 285-288. 
    
\item  Kariya, T., Giri, N., \& Perron, F. (1990).  Invariant estimation of a mean vector $\mu$ of $N(\mu,\Sigma)$ with
$\mu' \Sigma^{-1} \mu=1$, $\Sigma^{-1/2} \mu = C$ or $\Sigma= \delta^2 \mu' \mu I$.  {\it Journal of Multivariate Analysis},
{\bf 27}, 270-283.

\item Kokologiannaki, C.G. (2012).  Bounds for functions involving ratios of Bessel functions.  
{\it Journal of Mathematical Analysis and Applications}, {\bf 385}, 737-742.

\item  Kortbi, O. \& Marchand, \'E.  (2012).  Truncated linear estimation of a bounded multivariate normal mean,
{\it Journal of Statistical Planning and Inference}, http://dx.doi.org/10.1016/j.jspi.2012.03.022
 
%\item  Kortbi, O. (2011).  {\it Sur l'estimation d'une moyenne multivari\'ee sous sym\'etrie sph\'erique et sous contrainte.}  %Ph.D thesis, Universit\'e de Sherbrooke.
\item Kubokawa, T. (2005).  Estimation of a mean of a normal
distribution with a bounded coefficient of variation. {\it
Sankhy$\bar{a}$: The Indian Journal of Statistics}, {\bf 67}, 499-525.

%\item Fourdrinier, D. \& Wells, M.T. (1996). 
%Spherically symmetric Bayes estimators for a linear subspace of a normal law.
%{\it Bayesian Statistics 5},  Oxford University Press, New York, pp. 569-579
\item Kubokawa, T. (1994).  A unified approach to improving on equivariant estimators.
{\it Annals of Statistics}, {\bf 22}, 290-299.

\item  Joshi, C.M. \& Bissu, S.K. (1996).  Inequalities for some special functions.
{\it Journal of Computational and Applied Mathematics}, {\bf 69},  251-259.

\item  Marchand, \'{E}., Jafari Jozani, M. \& Tripathi, Y. M. (2012). On the inadmissibility of various estimators of
normal quantiles and on applications to two-sample problems with
additional information. {\it Contemporary
Developments in Bayesian analysis and Statistical Decision Theory: 
A Festschrift for William E. Strawderman}, Institute of
Mathematical Statistics Volume Series, {\bf 8}, 104-116.

\item  Marchand, \'{E}. \& Perron, F. (2001). Improving on the MLE of
a bounded normal mean. {\it Annals of Statistics}, {\bf 29}, 1078-1093.

\item  Marchand, \'E., and Strawderman, W.E. (2012).  
A unified minimax result for restricted parameter spaces.  
{\it Bernoulli}, {\bf 18}, 635-643.

\item  Marchand, \'{E}., and Strawderman, W.E. (2005).  On
improving on the minimum risk equivariant estimator of a location
parameter which is constrained to an interval or a half-interval.
{\it Annals of the Institute of Statistical Mathematics}, {\bf
57}, 129-143.

\item
Marchand, \'{E}. and Strawderman, W. E. (2004). Estimation in
restricted parameter spaces: A review. {\it Festschrift for Herman
Rubin}, IMS Lecture Notes-Monograph Series, {\bf 45}, pp. 21-44.

\item Moors, J.J.A. (1985).  Estimation in truncated parameter spaces.  Ph.D. thesis, Tilburg University.

\item  Perron, F. and Giri, N.C. (1990).  On the best equivariant estimator of a multivariate normal population.  
{\it Journal of Multivariate Analysis}, {\bf 32}, 1-16.

\item  Watson, G.S. (1983).  {\it Statistics on spheres}.  John Wiley, New York.

\end{description}

\end{document}